\documentclass[a4paper,reqno,9pt]{amsart}
\usepackage[margin=2.5cm,papersize={19.75cm,28cm}]{geometry}

\usepackage{amsmath}
  \usepackage{paralist}
  \usepackage{graphics} 
  \usepackage{epsfig} 
\usepackage{graphicx}  \usepackage{epstopdf}
 \usepackage[colorlinks=true]{hyperref}
\usepackage[all]{xy} 
\usepackage{amssymb} 
\hypersetup{urlcolor=blue, citecolor=red}
\usepackage{multicol}
\usepackage[active]{srcltx}
\usepackage{amsthm}
\usepackage{amsmath}
\usepackage{wasysym}
\newtheorem{theorem}{Theorem}[section]

\theoremstyle{definition}
\newtheorem{definition}[theorem]{Definition}
\newtheorem{remark}{Remark}

\newtheorem{example}{Example}

\newcommand{\ab}{\lbrack\!\lbrack \;, \; \rbrack\!\rbrack}

\def\lcf{\lbrack\! \lbrack}
\def\rcf{\rbrack\! \rbrack}
\newcommand{\lp}{\left(}
\newcommand{\rp}{\right)}
\newcommand{\lc}{\left\{}
\newcommand{\rc}{\right\}}

\newcommand{\R}{\mathbb{R}}      

\usepackage{amssymb}

\newcommand{\proa}{A^*G \mbox{$\;$}_{\tau^*} \kern-3pt\times_\alpha
G \mbox{$\;$}_\beta \kern-3pt\times_{\tau^*} A^*G}

\newcommand{\al}{\mathfrak{g}}

\title[Variational-geometric optimal control for nonholonomic systems]{A variational-geometric approach for the optimal control of nonholonomic systems}

\author[Leonardo Colombo]{}


\begin{document}
\maketitle

\centerline{\scshape Leonardo Colombo}
\medskip
{\footnotesize
 \centerline{Department of Mathematics, University of Michigan}
   \centerline{530 Church Street, 3828 East Hall}
   \centerline{ Ann Arbor, Michigan, 48109, USA}
   \centerline{ljcolomb@umich.edu}
} 

\bigskip




\begin{abstract}
Necessary conditions for existence of normal extremals in optimal control of systems subject to nonholonomic constraints are derived as solutions of a constrained second order variational problems. In this work, a geometric interpretation of the derivation is studied from the theory of Lie algebroids. We employ such a framework to describe the problem into a unifying formalism for normal extremals in optimal control of nonholonomic systems and including situations that have not been considered before in the literature from this perspective. We show that necessary conditions for existence of extremals  in the optimal control problem can be also determined by a Hamiltonian system on the cotangent bundle of a skew-symmetric algebroid.

\end{abstract}

\section{Introduction}
A nonholonomic system is a mechanical system subject to constraint functions which are, roughly speaking, functions on the velocities that are not derivable from position constraints.
They arise, for instance, in mechanical systems that have rolling or certain kinds
of sliding contact. There are multiple applications in the context of wheeled
motion, mobile robotics and robotic manipulation \cite{Bl, jorgebook, NF}. 

Optimal control problems of nonholonomic systems arise in many engineering applications, for instance, systems with wheels, such as maneuvers with cars or  bicycles, systems with blades or skates, and spherical robots \cite{blochcrouch3, hussein, koon}. We are mainly interested in the study of fully actuated systems, that is, when the number of control inputs is equal to the rank of the control distribution.

The goal of this note is to study, from a variational and geometric framework, necessary conditions for the existence of normal extremals in the optimal control of nonholonomic systems,  giving rise to a  unified framework that can include systems with phase space given by tangent bundles, Lie algebras, principal bundles, systems with symmetries as well as nonholonomic systems, instead of work case by case. To give such an approach we choice as theoretical framework of this work the theory of Lie algebroids \cite{CoLeMaMa}.

Our derivation of necessary conditions for normal extremals is determined by studying optimal control problems as constrained higher-order variational problems \cite{Bl, blochcrouch2, margarida}. These constrained higher-order variational problems are determined by minimizing the action associated to a $k^{th}$-order Lagrangian function defined on the $k^{th}$-order tangent bundle of a smooth manifold \cite{LR1},  giving rise to a $2k$-order system of ordinary differential equations and $k^{th}$-order constraints (subject to $2k$ boundary conditions). 

The results of this work employs the framework proposed in \cite{maria2} for kinematic mechanical control systems on skew-symmetric algebroids. This geometric structure allow us to describe in a unified and simple way the dynamics of non-standard (i.e., defined on tangent bundles) nonholonomic systems and the related optimal control problem. It is important to highlight that this
is not an arbitrary generalization since mechanics on algebroids is particularly relevant for the class of Lagrangian systems invariant under the action of a Lie group of symmetries including as a particular case nonholonomic dynamics and systems defined on Lie algebras and principal bundles \cite{CoLeMaMa,CoMa,januzleones,MdLMdDJCM}. 

The examples studied in this paper includes mechanical systems on Lie algebras, a situation that with our previous related work \cite{BlCoGuMdD} we can not study since we are restricted to work on the tangent bundles of the configuration manifold. Here, we avoid that obstacle in the phase space of the systems by considering the framework of (skew-symmetric) Lie algebroids. Therefore, the results of this note must be considered an extension of our previous work \cite{BlCoGuMdD}, presenting a general framework for optimal control of nonholonomic systems that allows to incorporate new situations than the ones studied previously in \cite{BlCoGuMdD}, also giving rise to a new contribution for applications of theories developed on Lie algebroids besides just contributing to our previous work \cite{BlCoGuMdD}. 

Necessary conditions studied in this notes are reduced equations, in the sense that the order of the equation we obtain for normal extremals using this approach must satisfy a first order differential equation on the constraint distribution, instead of a 4th order ordinary differential equation on the configuration manifold as is usual from the approach of higher-order variational calculus, and by using an admisible condition for the curves that satisfying the constraint, it is possible reconstruct solutions to the configuration manifold.

We also derive the corresponding Hamiltonian representation of optimal control problem when the system is regular. That framework permits to describe necessary conditions for regular extremal as solutions of a Hamiltonian systems defied on a symplectic manifold, the cotangent bundle of the nonholonomic distribution and one can then use standard methods for symplectic integration to integrate numerically the equations determining extremals for the control system if it is needed.  We show that the techniques of this work can be easily adapted for underactuated systems. The application of our techniques is tested in two nonholonomic systems on Lie algebras: The Suslov problem and the Chaplyging sleigh. 



\section{Lagrangian System of Mechanical Type}
Let $Q$ be a differentiable manifold of dimension $n$, the configuration space of a mechanical system,  with local coordinates $(q^A)$, $A=1,\ldots,n$, and  ${\mathcal G}$ be a Riemannian metric specifying the kinetic energy of the system. The metric is locally written as ${{\mathcal G}} = {{\mathcal G}}_{AB} dq^A\otimes dq^B$, where $\mathcal{G}_{AB}=\mathcal{G}(\partial/\partial q^A,
\partial/\partial q^B)$. Using the Riemannian metric it is possible to construct  the Levi-Civita connection on $Q$, $\nabla^{{\mathcal G}}:\mathfrak{X}(Q)\times\mathfrak{X}(Q)\to\mathfrak{X}(Q)$, where $\mathfrak{X}(Q)$ denotes the set of vector fields on $Q$, as the unique affine connection which is torsion-less and metric with respect to ${\mathcal G}$. It is determined by the standard formula \begin{align*}
2 {{\mathcal G}}(\nabla_{X}^{{\mathcal G}} Y, Z)=&
X({{\mathcal G}}(Y, Z)) + Y({{\mathcal G}}(X, Z))- Z({{\mathcal G}}(X, Y)) \\
& +  {{\mathcal G}}(X, [Z, Y]) + {{\mathcal G}}(Y, [Z,
X]) - {{\mathcal G}}(Z, [Y, X])
\end{align*} for all $X, Y, Z \in\mathfrak{X}(Q)$ \cite{jorgebook}. Alternatively, $\nabla^{{\mathcal G}}$ is determined by the properties of symmetry and metricity for the connection: $[X,Y]=\nabla^{{\mathcal G}}_X Y-\nabla^{{\mathcal G}}_Y X$ and
$X({{\mathcal G}}(Y,Z))={{\mathcal G}}(\nabla^{{\mathcal G}}_X Y, Z)+{{\mathcal G}}(Y, \nabla^{{\mathcal G}}_X Z)$ respectively.


Fixed a potential  function $V: Q\rightarrow \R$,  the mechanical system
is defined by the \textit{mechanical Lagrangian} $L: TQ\rightarrow \mathbb{R}$,
\begin{equation}\label{lagrangian}
L(v_q)=\frac{1}{2}{{\mathcal G}}(v_q,v_q)-V(q),
\end{equation} where $v_q\in T_qQ$ and the solutions of the variational problem
determined by $L: TQ\rightarrow \R$ are curves $c: I\subset \R\to Q$ such that
\begin{equation}\label{equations}
\nabla^{{\mathcal G}}_{\dot{c}(t)}\dot{c}(t) + \hbox{grad}_{{\mathcal
G}}V(c(t))=0.
\end{equation}
Here, $\hbox{grad}_{{\mathcal G}}V$ is the vector field on $Q$ characterized by $
{{\mathcal G}}(\hbox{grad}_{{\mathcal G}}V, X) =
X(V)$, for  every  $X \in \mathfrak{X}(Q)$. If $V=0$ then $c(t)$ is the solution of the geodesic equations.

 In local coordinates, equations (\ref{equations}) are given by
\begin{equation}
\ddot{q}^{C}=-\Gamma^C_{AB}(q(t))\dot{q}^{A}\dot{q}^{B}-{{\mathcal G}}^{AB}\frac{\partial V}{\partial q^C}\label{equa1}\;, 
\end{equation} where $({{\mathcal G}}^{AB})$ are the entries of the inverse matrix of $({\mathcal G}_{AB})$  and where $\Gamma_{AB}^C$ are the Christoffel symbols associated with the Levi-Civita connection and computed from the formula
$\displaystyle{\nabla^{{\mathcal G}}_{\frac{\partial}{\partial q^A}}{\frac{\partial}{\partial q^B}}=\Gamma^C_{AB} \frac{\partial}{\partial q^C}}$.
If $V=0$ the equation reduces to the local description of the geodesic equations
 $\ddot{q}^C=-\Gamma^C_{AB}\dot{q}^A\dot{q}^B$.

\section{Geometry of nonholonomic mechanical systems}
In mechanics usually appear two type of constraints: holonomic and nonholonomic constraints. A \textit{holonomic} constraint restricts the dynamics only in terms of position, or in other words, it tells where the dynamics should be, while a \textit{nonholonomic} constraint does it in terms of velocity only, or it tells in which direction the dynamics should go. Typical examples of nonholonomic constraints are those imposed by rolling and sliding of the mechanical systems. Such systems often arise in engineering problems, e.g., systems with wheels like cars and bicycles and those with sliding parts like sleighs \cite{Bl, NF}. Next, we will study the underlying geometry of the dynamics described by nonholonomic systems.

\subsection{Nonholonomic mechanical systems on tangent bundles}

A \textit{nonholonomic system} is a mechanical system with external constraints on the velocities. 
We only consider \textit{linear velocity constraints}, since this is the case in most examples. 

Linear velocity constraints are constraints that are specified by a regular $C^{\infty}$-distribution ${\mathcal D}$  on the configuration manifold $Q$, or equivalently, by a vector subbundle $\tau_{\mathcal D}: {\mathcal D}\rightarrow Q$ of $TQ$ where the inclusion is denoted by $i_{\mathcal D}: {\mathcal D}\hookrightarrow TQ$. Therefore, we will say that a curve $\gamma: I\subseteq {\mathbb R}\rightarrow Q$ \textit{satisfies the constraints} given by ${\mathcal D}$ if 
$\dot{\gamma}(t)\in {\mathcal D}_{\gamma(t)}$ for all $t\in I$.

We say that ${\mathcal D}$ is holonomic if ${\mathcal D}$ is \textit{integrable} or involutive, that is, for any vector fields  $X, Y\in {\mathfrak X}(Q)$ taking values on ${\mathcal D}$, it holds that the vector field $[X, Y]$ also takes values on ${\mathcal D}$. A regular linear velocity constraint submanifold ${\mathcal D}$ is nonholonomic if  it is not holonomic. 
Observe that in the case of holonomic constraints all the curves through a point $q\in Q$ satisfying the constraints must lie on the maximal integral manifold for ${\mathcal D}$ through $q$. 

Let $\dim Q=n$. Locally, if $(q^A)$, $1\leq A\leq n$ are coordinates on $Q$ and 
$(q^A, \dot{q}^A)$ are the induced coordinates on $TQ$,  the linear constraints are written as 
\begin{equation*}\label{LC}
\mu^{\alpha}_{A}\lp q\rp\dot q^{A}=0,\hspace{2mm} m+1\leq \alpha\leq n\, ,
\end{equation*}
where $\mbox{rank}\lp\mathcal{D}\rp=m\leq n$. The annihilator $\mathcal{D}^{\circ}$ is locally given by
\[
\mathcal{D}^{\circ}=\mbox{span}\lc\mu^{\alpha}=\mu_{A}^{\alpha}(q)\,dq^{A};\hspace{1mm} m+1\leq \alpha\leq n\rc\, ,
\]
where the 1-forms $\mu^{\alpha}$ are independent.
Equivalently, we can find independent vector fields $\{X_a\}$, $1\leq a\leq m$ such that ${\mathcal D}_q=\hbox{span}\{X_a\}\,$. Observe that $\mu^{\alpha}(X_a)=0$, for all $m+1\leq \alpha\leq n$ and $1\leq a\leq m$.

Now we restrict ourselves to nonholonomic mechanical systems where the Lagrangian is of \textit{mechanical type}, that is, a Lagrangian systems $L:TQ\to\R$ defined by \[
L(v_q)=\frac{1}{2}\mathcal{G}(v_q, v_q) - V(q),
\]
with $v_q\in T_qQ$, $\mathcal{G}$ denotes a Riemannian metric on the configuration
space $Q$ representing the kinetic energy of the systems and
$V:Q\to\R$ is a potential function, as in Section $2$.


A \textit{nonholonomic mechanical system} on a smooth manifold $Q$ is given
by the triple $(\mathcal{G}, V, \mathcal{D})$ where  $\mathcal{G}$ and $V$ as before and $\mathcal{D}$ a non-integrable regular distribution on $Q$.

Denoting by $\mathfrak{X}(\mathcal{D})$ the set of vector fields taking values on $\mathcal{D}$ and $\mathfrak{X}(Q)$ the one taking values on $TQ$, if $X, Y\in\mathfrak{X}(\mathcal{D})$ then
$[X,Y]$ denotes the standard Lie bracket of vector fields. Given $X,Y\in\mathfrak{X}(\mathcal{D})$ that is,
$X(x)\in\mathcal{D}_{x}$ and $Y(x)\in\mathcal{D}_{x}$ for all $x\in
Q,$ then it may happen that $[X,Y]\notin\mathfrak{X}(\mathcal{D})$
since $\mathcal{D}$ is nonintegrable. 

In order to obtain a bracket definition for vector field taking values on $\mathcal{D}$ (and therefore satisfying the constraints) one may uses the Riemannian metric $\mathcal{G}$ to define two
complementary orthogonal projectors ${\mathcal P}\colon TQ\to {\mathcal D}$ and ${\mathcal Q}\colon TQ\to {\mathcal
D}^{\perp},$ with respect to the tangent bundle orthogonal decomposition $\mathcal{D}\oplus\mathcal{D}^{\perp}=TQ$. Therefore, given $X,Y\in\mathfrak{X}(\mathcal{D})$ we define a new bracket, 
$\lcf\cdot,\cdot\rcf:\mathfrak{X}(\mathcal{D})\times\mathfrak{X}(\mathcal{D})\rightarrow\mathfrak{X}(\mathcal{D})$
as $\lcf X,Y\rcf:=\mathcal{P}[X,Y]$. This Lie bracket verifies the usual properties of a Lie bracket, except, in particular, the Jacobi identity \cite{CoMa}.

\begin{definition}
A curve $\gamma:I\subset\R\to\mathcal{D}$ is \textit{admissible} if
$\gamma(t)=\dot{\sigma}(t)$, where $\tau_{\mathcal{D}}\circ\gamma=\sigma$.
\end{definition}

Given local coordinates on $Q,$ $(q^{i})$ with $i=1,\ldots,n;$ and
$\{e_{A}\}$ a basis of vecotr fields on $\mathfrak{X}(\mathcal{D})$, with $A=1,\ldots,n-m$, such that
$\displaystyle{e_{A}=\rho_{A}^{i}(q)\frac{\partial}{\partial
q^{i}}}$ we introduce induced coordinates $(q^{i},y^{A})$ on
$\mathcal{D}$, where, if $e\in\mathcal{D}_{x}$ then
$e=y^{A}e_{A}(x).$ Therefore, $\gamma(t)=(q^{i}(t),y^{A}(t))$ is
admissible if
$\dot{q}^{i}(t)=\rho_{A}^{i}(q(t))y^{A}(t).$

Consider the restriction of the Riemannian metric $\mathcal{G}$ to
the distribution $\mathcal{D}$, denoted by $\mathcal{G}^{\mathcal{D}}:\mathcal{D}\times_{Q}\mathcal{D}\to\R$. 
The \textit{Levi-Civita connection}
$\displaystyle{\nabla^{\mathcal{G}^{\mathcal{D}}}:\mathfrak{X}(\mathcal{D})\times\mathfrak{X}(\mathcal{D})\to\mathfrak{X}(\mathcal{D})}$
is determined by the symmetry property
$\lcf X,Y\rcf=\nabla_{X}^{\mathcal{G}^{\mathcal{D}}}Y-\nabla_{Y}^{\mathcal{G}^{\mathcal{D}}}X,$ and the metricity $X(\mathcal{G}^{\mathcal{D}}(Y,Z))=\mathcal{G}^{\mathcal{D}}(\nabla_{X}^{\mathcal{G}^{\mathcal{D}}}Y,Z)+\mathcal{G}^{\mathcal{D}}(Y,\nabla_{X}^{\mathcal{G}^{\mathcal{D}}}Z).$

\begin{definition}\cite{maria2}
Consider the restricted Lagrangian function
$\ell:\mathcal{D}\rightarrow\mathbb{R},$
$$\ell(v)=\frac{1}{2}\mathcal{G}^{\mathcal{D}}(v,v)-V(\tau_{D}(v)),\hbox{ with }  v\in\mathcal{D}.$$

A \textit{solution of the nonholonomic problem} is an admissible
curve $\gamma:I\rightarrow\mathcal{D}$ such that
$$\nabla_{\gamma(t)}^{\mathcal{G}^{\mathcal{D}}}\gamma(t)+grad_{\mathcal{G}^{\mathcal{D}}}V(\tau_{\mathcal{D}}(\gamma(t)))=0,$$
where the vector field $grad_{{\mathcal G}^{\mathcal{D}}}V\in\mathfrak{X}(\mathcal{D})$ is characterized by $ {{\mathcal G}^{\mathcal{D}}}(grad_{{\mathcal
G}^{\mathcal{D}}}V, X) = X(V)$ for  every  $X \in
\mathfrak{X}(\mathcal{D})$.
\end{definition} 

Locally, admisible solutions for the nonholonomic problem are determined by
$$
\dot{q}^{i}=\rho_{A}^{i}(q)y^{A},\qquad
\dot{y}^{C}=-\Gamma_{AB}^{C}y^{A}y^{B}-(\mathcal{G}^{\mathcal{D}})^{CB}\rho_{B}^{i}\frac{\partial V}{\partial q^{i}},$$ where $(\mathcal{G}^{\mathcal{D}})^{AB}$ denotes the coefficients of the inverse matrix of $(\mathcal{G}^{\mathcal{D}})_{AB}$ where $\mathcal{G}^{\mathcal{D}}(e_{A},e_{B})=(\mathcal{G}^{\mathcal{D}})_{AB}$, and the \textit{Christoffel
symbols} $\Gamma_{BC}^{A}$ of the connection
$\nabla^{\mathcal{G}^{\mathcal{D}}}$ can be determined by $\displaystyle{\nabla_{e_{B}}^{\mathcal{G}^{\mathcal{D}}}e_{C}=\Gamma_{BC}^{A}(q)e_{A}}$. 
\subsection{Nonholonomic systems on Lie algebroids}
Instead of work on $TQ$ we can consider an arbitrary Lie algebroid $E$ \cite{januzleones}.
The projection onto a real vector subbundle $\mathcal{D}\subset E$ of the Lie bracket determined by the Lie algebroid structure of $E$ gives rise to a skew-symmetric algebroid structure on $\mathcal{D}$ \cite{MdLMdDJCM}. This approach permits to include in the analysis systems with Lie algebras and principal bundles as phase space, situations that does not allow to include our previous work on $TQ$ \cite{CoBl}.

It is known that this geometric
structure covers many interesting cases in mechanics, as for
instance, nonholonomic mechanics for systems defined on Lie algebra, principal bundles, and reduced systems \cite{CoLeMaMa}, \cite{januzleones}, \cite{MdLMdDJCM}. Similarly to
the intrinsic definition of the Euler-Lagrange equations for a
Lagrangian function $L: TQ\to \R$ obtained by the canonical
structures on it (standard Lie bracket, exterior differential...),
it is possible to derive the dynamics of the system from a Lagrangian
$L: \mathcal{D}\to \R$ using the differential geometric structures naturally
induced by the skew-symmetric algebroid structure.
This generalization is useful in applications and clarifies
the dynamics of systems with nonholonomic constraints. 

Nevertheless, in this note, we prefer to use as a starting point the Lie algebroid structure on $E$ and derive the skew-symmetric algebroid structure instead of start with a skew-symmetric algebroid structure as in \cite{maria2}. Both formalisms are dynamically equivalent.

A Lie algebroid $E$ of rank $n$ over a
manifold $Q$ of dimension $m$, is a real vector bundle $E$ with projection
$\tau_E:E\rightarrow Q$ together with a Lie bracket $\ab_{E}$ on
$\Gamma(\tau_E)$, the set of sections of $\tau_{E}:E\to Q$,  and
a fiber map $\rho_E:E\to TQ$ called \textit{anchor map}. We will denote the Lie algebroid $E$ by the triple $(E,\ab_E,\rho_E).$ It would be helpful for readers without previous background on Lie algebroids think sections of $\tau_E$, as vector fields on $Q$, and the sections of the dual bundle
$\tau_{E^{*}}:E^{*}\to Q$, like 1-forms on $Q$.
%

\begin{definition}
A \textit{nonholonomic system} on a Lie algebroid
$(E,\rho_{E},\lcf\cdot,\cdot,\rcf_{E})$ over a manifold $Q$ with
bundle projection $\tau_{E}:E\to Q$ is a triple $(\mathcal{D},\mathcal{G},V)$ determined by the following
three data: a real vector subbundle $\mathcal{D}$ of $E$, a nondegenerate bundle metric $\mathcal{G}$ on $E$, $\mathcal{G}:E\times_{Q}E\to\R$, and a smooth function $V:Q\to\R$.

\end{definition}

Using the bundle metric it is possible to construct two complementary
projectors, $\mathcal{P}:E\to\mathcal{D}$ and 
$\mathcal{Q}:E\to\mathcal{D}^{\perp}$, with respect to the orthogonal decomposition
$E=\mathcal{D}\oplus\mathcal{D}^{\perp}$, and projecting the Lie bracket on $\Gamma(\tau_{E})$ to $\mathcal{D}$, we obtain a new Lie bracket over
$\Gamma(\tau_{\mathcal{D}})$ as \[\lcf
X,Y\rcf_{\mathcal{D}}:=\mathcal{P}\lcf
i_{\mathcal{D}}(X),i_{\mathcal{D}}(Y)\rcf_{E},
\] where $X,Y\in\Gamma(\tau_{\mathcal{D}}),$ $\tau_{\mathcal{D}}:\mathcal{D}\to Q$ is the restriction of $\tau_E$ to $\mathcal{D}$ and $i_{D}:\mathcal{D}\to E$ is
the inclusion of the subbundle $\mathcal{D}$ on $E$.

Denoting local coordinates on $Q$ by $(q^i)$ and $\{e_{A}\}$ be
a local basis of the space of sections $\Gamma(\tau_\mathcal{D})$,
then
\begin{equation*} \lcf e_{A}, e_B\rcf_\mathcal{D} =
{\mathcal C}^C_{AB} e_{C}, \ \ \ \rho_\mathcal{D} (e_{A})=(\rho_\mathcal{D})_{A}^i \frac{\partial}{\partial q^i},  \label{coeff-estruct}
\end{equation*} where $\rho_{\mathcal{D}}:\mathcal{D}\to TQ$ is the restriction of
$\rho_{E}$ to $\mathcal{D}$ satisfying $\rho_{\mathcal{D}}(X)=i_{\mathcal{D}}(X)$ for $X\in\Gamma(\tau_{\mathcal{D}})$. 

The triple $(\mathcal{D}, \lcf
\cdot,\cdot\rcf_{\mathcal{D}}, \rho_{\mathcal{D}})$ is know as \textit{skew-symmetric Lie algebroid} \cite{januzleones, MdLMdDJCM} and the functions ${\mathcal C}^C_{AB}, (\rho_D)_A^i\in C^{\infty}(Q)$
are called the \emph{local structure functions} of $(\mathcal{D}, \lcf
\cdot,\cdot\rcf_{\mathcal{D}}, \rho_{\mathcal{D}})$.

A $\rho_{\mathcal{D}}$\emph{-admissible curve} is a curve $\gamma: I\subseteq
\mathbb{R} \longrightarrow \mathcal{D}$  such that
\[
\frac{d}{dt}(\tau_{\mathcal{D}}\circ \gamma)=\rho_{\mathcal{D}}(\gamma(t))\;
.
\] Locally, if we take local coordinates $(q^i)$ on $Q$ and a basis of sections $\{e_A\}$ of $\tau_{\mathcal{D}}$, then we have the
corresponding induced coordinates $(q^i,y^A)$ on $\mathcal{D}$,
where $y^A(a)$ is the $A$-th coordinate of $a\in \mathcal{D}$ in the
given basis. Therefore, $\gamma(t)=(q^{i}(t),y^{A}(t))$ is
$\rho_{\mathcal{D}}$-admissible if
$$\dot{q}^{i}=(\rho_{\mathcal{D}})_{A}^{i}y^{A}.$$ Moreover, given $X\in \Gamma(\tau_{\mathcal{D}})$, the \emph{integral curves}
of the section $X$ are those curves $\sigma: I\subseteq
\R\rightarrow Q$ such that satisfy $\dot{\sigma}=\rho_{\mathcal{D}}(X)\circ \sigma$. That is, they are the integral curves of the associated vector field
$\rho_{\mathcal{D}}(X)\in {\mathfrak X}(Q)$. 

If $\sigma$ is an integral curve of
$X$, then $X\circ \sigma$ is a $\rho_{\mathcal{D}}$-admissible curve. Locally,
the integral curves are characterized as the solutions of the
system of equations $\dot{q}^i=(\rho_{\mathcal{D}})^i_A X^A(q)$, where $X=X^Ae_A$ \cite{maria2}.

%
The bundle metric restricted to the vector subbundle $\mathcal{D}$, denoted by $\mathcal{G}^{\mathcal{D}}:\mathcal{D}\times_Q\mathcal{D}\to\R$, and locally determined by $\mathcal{G}^{\mathcal{D}}=(\mathcal{G}^{\mathcal{D}})^{AB}e^{A}\otimes e^{B}$,permits to construct a unique torsion-less
connection $\nabla^{{\mathcal G}^{\mathcal{D}}}$ on $\mathcal{D}$. The \textit{Levi-Civita connection} $\nabla^{{\mathcal
G}^\mathcal{D}}: \Gamma(\tau_\mathcal{D})\times
\Gamma(\tau_\mathcal{D})\to \Gamma(\tau_\mathcal{D})$ associated to
the bundle metric ${{\mathcal G}^{\mathcal{D}}}$ is defined  by the
formula \begin{align*}
2 {{\mathcal G}^{\mathcal{D}}}(\nabla_{X}^{{\mathcal G}^{\mathcal{D}}} Y, Z) = &
\rho_{\mathcal{D}}(X)({{\mathcal G}^{\mathcal{D}}}(Y, Z)) + \rho_{\mathcal{D}}(Y)({{\mathcal G}^{\mathcal{D}}}(X, Z))- \rho_{\mathcal{D}}(Z)({{\mathcal G}^{\mathcal{D}}}(X, Y)) \\
& +  {{\mathcal G}^{\mathcal{D}}}(X, \lcf Z, Y\rcf_{\mathcal{D}})
+ {{\mathcal G}^{\mathcal{D}}}(Y,\lcf Z,
X\rcf_{\mathcal{D}}) - {{\mathcal G}^{\mathcal{D}}}(Z, \lcf Y, X\rcf_{\mathcal{D}})
\end{align*} for $X, Y, Z \in \Gamma(\tau_\mathcal{D})$. Alternatively, $\nabla^{{\mathcal G}^{\mathcal{D}}}$ is determined
by the properties of symmetry $ \lcf X,Y\rcf_{\mathcal{D}}= \nabla^{{\mathcal G}^{\mathcal{D}}}_X Y-\nabla^{{\mathcal G}^{\mathcal{D}}}_Y X$ and metricity
$\rho_\mathcal{D}(X)({{\mathcal G}^{\mathcal{D}}}(Y,Z))={{\mathcal G}^{\mathcal{D}}}(\nabla^{{\mathcal G}^{\mathcal{D}}}_X Y, Z)+{{\mathcal G}^{\mathcal{D}}}(Y, \nabla^{{\mathcal G}^{\mathcal{D}}}_X Z)$. Usually, the Levi-Civita connection $\nabla^{\mathcal{G}^{\mathcal{D}}}$
coincides with the constrained connection $\nabla^{\mathcal{D}}_{X}Y:=\mathcal{P}(\nabla^{\mathcal{G}}_{X}Y)$ defined for instance in \cite{CoMa}, if $\nabla^{\mathcal{D}}$ is restricted to $\Gamma(\tau_{\mathcal{D}})$.

As when we work in tangent bundles, it is possible to determine the \textit{Christoffel
symbols} associated with the connection $\nabla^{{\mathcal
G}^{\mathcal{D}}}$ by $\nabla^{{\mathcal G}^{\mathcal{D}}}_{e_B}{e_C}=\Gamma^A_{BC}e_A$. Note that the coefficients $\Gamma_{AB}^{C}$ of the connection  $\nabla^{{\mathcal
G}^{\mathcal{D}}}$ are (see \cite{CoLeMaMa,CoMa} for details) \begin{equation}\label{relation}\Gamma_{AB}^{C}=\frac{1}{2}(\mathcal{C}_{CA}^{B}+\mathcal{C}_{CB}^{A}+\mathcal{C}_{AB}^{C}).\end{equation}

\begin{definition}
A \textit{solution of the nonholonomic problem} is a
$\rho_{\mathcal{D}}$-admissible curve
$\gamma:I\subset\R\rightarrow\mathcal{D}$ such that
$$\nabla_{\gamma(t)}^{\mathcal{G}^{\mathcal{D}}}\gamma(t)+grad_{\mathcal{G}^{\mathcal{D}}}V(\tau_{\mathcal{D}}(\gamma(t)))=0.$$
\end{definition}

Here, $grad_{{\mathcal G}^{\mathcal{D}}}V$ is a section of
$\tau_{\mathcal{D}}:\mathcal{D}\to Q$ characterized by \[ {{\mathcal
G}^{\mathcal{D}}}(grad_{{\mathcal G}^{\mathcal{D}}}V, X) = \rho_{\mathcal{D}}(X)(V), \;
\; \mbox{ for  every } X \in \Gamma(\tau_{\mathcal{D}}).
\] Locally, solution must satisfy

\begin{equation}\label{equationlocal}\dot{q}^i=(\rho_{\mathcal{D}})^i_A y^A,\quad
\dot{y}^C=-\Gamma^C_{AB}y^Ay^B-({{\mathcal
G}^{\mathcal{D}}})^{CB}(\rho_{\mathcal{D}})^i_B\frac{\partial V}{\partial q^i}.\;\end{equation}

 \begin{example}[Euler-Poincar\'e-Suslov equations on $\mathfrak{so}(3)$]
 {\rm
As an example we study nonholonomic systems defined on a finite dimension real Lie
algebra ${\mathfrak g}$. It is well know that ${\mathfrak g}$ is a Lie
algebroid over a single point where the anchor map is $\rho= 0$ and the bracket is determined by the Lie algebra structure of $\al$ \cite{CoLeMaMa}. 

Consider a nonholonomic Lagrangian system on ${\mathfrak g}$, determined by
$L: {\mathfrak g}\to \R$, a metric (kinetic) Lagrangian function defined by
$L(\xi)=\frac{1}{2}\langle \mathbb{I} \xi, \xi\rangle$, with $\mathbb{I}: {\mathfrak
g}\rightarrow {\mathfrak g}^*$ a symmetric positive definite
inertia operator and ${\mathfrak D}$  a linear vector subspace of
${\mathfrak g}$. The orthogonal decomposition ${\mathfrak g}:={\mathfrak D}\oplus {\mathfrak D}^{\perp}$, with ${\mathfrak D}^{\perp}=\{\eta\in {\mathfrak g}\, |\, \langle
\mathbb{I} \eta, \xi\rangle=0 \; \forall \xi \in {\mathfrak D}\}$ permits to define the
associated orthogonal projector
$\mathcal{P}:\mathfrak{g}\to\mathfrak{D}$. The
bracket on $\mathfrak{D}$ is determined by $\lcf \cdot,\cdot\rcf_{\mathcal{D}}=\mathcal{P}[\cdot,\cdot].$ By considering an 
adapted basis ${\mathfrak D}=\hbox{span } \{e_A\}$ it is possible to induce coordinates $(y^{A})$ on $\mathfrak{D}$ and determine a restricted Lagrangian $\ell:\mathfrak{D}\to\mathbb{R}$. The
\emph{Euler-Poincar\'e-Suslov equations} \cite{Bl} for $\ell: {\mathfrak D}\to\mathbb{R}$
are
\[
\dot{y}^C=-\Gamma^C_{AB}y^Ay^B. \;
\]

If $\mathfrak{g}=\mathfrak{so}(3)$, the Lie algebra  skew-symmetric $3\times 3$ matrices, by considering the basis $\{e_1,e_2,e_3\}$ of $\mathfrak{so}(3)\simeq\mathbb{R}^{3}$ we can induce adapted coordinates $\xi=(\xi^1,\xi^2,\xi^3)$ on $\mathfrak{D}$ where an element $\xi\in\mathfrak{so}(3)$ is given by $\xi=\xi^1e_1+\xi^2e_2+\xi^3e_3$. The inertia tensor $\mathbb{I}$ is given by $$\mathbb{I}=\left(\begin{array}{ccc} I_{11}&0&I_{13}\\ 0&I_{22}&I_{23}\\ I_{13}&I_{23}&I_{33}
\end{array} \right),$$ and the Lagrangian $L:\mathfrak{so}(3)\to\mathbb{R}$ by $$L(\xi)=\frac{1}{2}(I_{11}(\xi^1)^{2}+I_{22}(\xi^2)^2+I_{33}(\xi^3)^2+2\xi^1\xi^3I_{13}+2\xi^{2}\xi^{3}I_{23}).$$

If the dynamics is subject to the linear nonholonomic constraint $a\xi^3=0$ where $a\in\mathfrak{so}(3)\simeq\mathbb{R}^{3}$ one  can choose the basis  $\{e_1,e_2,e_3\}$ such that $a=e_3$, hence,  the constraint distribution is given by the linear subspace $\mathcal{D}=\{\xi\in\mathfrak{so}(3)\mid \xi^{3}=0\}$. 

Instead of $\{e_1, e_2, e_3\}$ we take the basis of
$\mathfrak{so}(3)$ adapted to the orthogonal decomposition ${\mathcal D}\oplus
{\mathcal D}^\perp$;  determined by $\{X=(1,0,0), Y=(0,1,0), Z=I_{22}I_{13},I_{11}I_{23},-I_{11}I_{22})\}$ where $\mathcal{D}=\hbox{span}\{X, Y\}$ and $\mathcal{D}^{\perp}=\hbox{span}\{Z\}$. We take adapted coordinates $(y_1,y_2)$ on $\mathcal{D}$ relative to the basis given by $\{X,Y\}$ in such a way an element $y\in\mathcal{D}$ can be written as $y=y_1X+y_2Y$.  in this sense, one can obtain the restricted Lagrangian $\ell:\mathcal{D}\to\mathbb{R}$ given by $$\ell(y_1,y_2)=\frac{1}{2}\left(I_{11}y_1^{2}+I_{22}y_{2}^2\right).$$ 
The unique non-vanisihng structure constants of the projected bracket are $\displaystyle{\mathcal{C}_{12}^{1}=\frac{I_{13}}{I_{11}}}$ and $\displaystyle{\mathcal{C}_{12}^{2}=\frac{I_{23}}{I_{22}}}$. Using the relation \eqref{relation} between the Christoffel symbols and the constant of structures, we can derive the Euler-Poincar\'e-Suslov equations on $\mathfrak{so}(3)$ $$\dot{y}_1=-\left(\frac{I_{13}}{I_{11}}y_1+\frac{I_{23}}{I_{22}}y_2\right)y_2,\quad \dot{y}_2=\left(\frac{I_{13}}{I_{11}}y_1+\frac{I_{23}}{I_{22}}y_2\right)y_1.$$

 \hfill$\diamond$}
 \end{example}
%
\section{Optimal control of nonholonomic mechanical systems on  Lie
algebroids}\
In this section we study necessary conditions of optimality for the class of fully actuated nonholonomic systems from a variational framework and we show that under some regularity conditions it is possible to transform the optimal control problem into a Hamiltonian system on $T^{*}\mathcal{D}$. We derive the Hamiltonian dynamics determining necessary conditions for existence of extremals in the optimal control problem. The Hamiltonian point of view can be used to derive symplectic integrators to study the qualitative behavior of solutions. 
\subsection{Optimal control of fully actuated mechanical systems}
The purpose of this section is to study optimal control problems for a nonholonomic mechanical systems. We will restric ourselves to the case when the dimension of the input or control distribution is equal to the rank of $\mathcal{D}$.  If the rank of $\mathcal{D}$ is equal to the dimension of the control distribution, the system will be called a \textit{fully actuated nonholonomic system}. Also we shall assume that all the mechanical  control systems in this work are controllable \cite{CoMa}.

Let $(\mathcal{D}, \lcf
\cdot,\cdot\rcf_{\mathcal{D}}, \rho_{\mathcal{D}})$ be a skew-symmetric Lie algebroid over a manifold $Q$ and assume that the nonholonomic system determined by $({\mathcal
G}^{\mathcal{D}}, L, {\mathcal  D})$ also contains some input
section $Y_{1},\ldots,Y_{m}$ with $m=\hbox{ rank }\mathcal{D}$. Therefore the \textit{control distribution} is
given by the vector subbundle $\mathcal{D}_{(c)}:=span\{Y_{A}\},$ where
$Y_{A}\in\Gamma(\tau_{\mathcal{D}})$. We will denote 
 $\{Y^{A}\}$ its dual basis, a basis of
$\Gamma(\tau_{\mathcal{D}^{*}})$. The \textit{equations of motion for a nonholonomic system with input
sections} are 
\begin{equation}\label{noholonoma-1}
\nabla^{{\mathcal G}^{\mathcal{D}}}_{\gamma(t)}\gamma(t) +grad_{\mathcal{G}^{\mathcal{D}}}V(\tau_{\mathcal{D}}(\gamma(t)))\in {\mathcal D}_{(c)}(\gamma(t)), \quad \forall \; t\in
I \subseteq \mathbb{R},
\end{equation}
where $\gamma: I\subset\R\to \mathcal{D}$ is a
$\rho_{\mathcal{D}}$-admissible curve \cite{maria2}.

In terms of control inputs, Equation (\ref{noholonoma-1}) can be
rewritten as
\begin{equation}\label{noholonoma-1-1}
\nabla^{{\mathcal G}^D}_{\gamma(t)}\gamma(t)+grad_{\mathcal{G}^{\mathcal{D}}}V(\tau_{\mathcal{D}}(\gamma(t)))
=\sum_{A=1}^m u^A(t) Y_A(\tau_D(\gamma(t)))
\end{equation}
for $u: U \to \R^m$, the control inputs, with $U$ is an open subset of $\mathbb{R}$ including the origin. Equivalently, solutions of the fully actuated 
nonholonomic problem are characterized by the admissible curves
which solve

\[
\Big{\langle}\nabla_{\gamma(t)}^{\mathcal{G}^{\mathcal{D}}}\gamma(t)+grad_{\mathcal{G}^{\mathcal{D}}}V(\tau_{\mathcal{D}}(\gamma(t))),Y^{A}(\tau_{\mathcal{D}}(\gamma(t))) \Big{\rangle}=u^{A}(t)\]
Locally, solution must satisfy $$\dot{q}^i=(\rho_{\mathcal{D}})^i_A y^A,\quad
\dot{y}^C+\Gamma^C_{AB}y^Ay^B+({{\mathcal
G}^{\mathcal{D}}})^{CB}(\rho_{\mathcal{D}})^i_B\frac{\partial V}{\partial q^i}=u^{C}\; .$$
  
\begin{definition}
The triple $({{\mathcal G}^{\mathcal{D}}}, L,
{\mathcal D}_{(c)})$ is called a \textit{fully actuated nonholonomic
mechanical control system on the skew-symmetric Lie algebroid} $(\mathcal{D}, \lcf
\cdot,\cdot\rcf_{\mathcal{D}}, \rho_{\mathcal{D}})$.
\end{definition}

Given a cost function
\begin{eqnarray*}
C:\mathcal{D}\times U&\longrightarrow&\mathbb{R}\\
(q^{i},y^{A},u^{a})&\longmapsto& C(q^{i},y^{A},u^{a})
\end{eqnarray*} the \textit{optimal control problem} consists on finding a $\rho_{\mathcal{D}}$-admissible curve
$\gamma:I\rightarrow\mathcal{D}$ solution of the fully actuated 
nonholonomic problem 
given boundary conditions on $\mathcal{D}$ and minimizing the cost
functional
$$\mathcal{J}(\gamma(t),u(t)):=\int_{0}^{T}C(\gamma(t),u(t))dt.$$
In order to find necessary conditions for optimal extremals, consider the subbundle $\mathcal{D}^{(2)}$ of $T\mathcal{D}$
\begin{equation*}
\mathcal{D}^{(2)}:=\{v\in T\mathcal{D}\mid v=\dot{\gamma}(0)\hbox{
where } \gamma:I\rightarrow\mathcal{D} \hbox{ is admissible}\}.
\end{equation*}
Locally $\mathcal{D}^{(2)}$ is described by the vanishing
of the constraints $\dot{q}^{i}-(\rho_{\mathcal{D}})_{A}^{i}y^{A}=0 \hbox{ on
}T\mathcal{D},$ where local coordinates on $T\mathcal{D}$ are
$(q^{i},y^{A},\dot{q}^{i},\dot{y}^{A})$ and coordintes on $\mathcal{D}^{(2)}$ are determined by $(q^i,y^A,\dot{y}^{A})$ where the inclusion from $\mathcal{D}^{(2)}$ to $T\mathcal{D}$, denoted by $i_{\mathcal{D}^{(2)}}:\mathcal{D}^{(2)}\to T\mathcal{D}$ is given by $$i_{\mathcal{D}^{(2)}}(q^{i},y^{A},\dot{y}^{A})=(q^i,y^{A},(\rho_\mathcal{D})_{A}^{i}y^{A},\dot{y}^{A}).$$

Solving the fully actuated nonholonomic control problem is
equivalent to solving a constrained second-order variational problem, determined by the Lagrangian
$\mathcal{L}:\mathcal{D}^{(2)}\rightarrow\mathbb{R}$ given, in the
selected coordinates, by 
$$\mathcal{L}(q^{i},y^{A},\dot{y}^{A})=C\left(q^{i},
y^{A},\dot{y}^{C}+\Gamma_{AB}^{C}y^{A}y^{B}+(\mathcal{G}^{\mathcal{D}})^{CB}\rho_{B}^{i}\frac{\partial
L}{\partial q^{i}}\right),$$ and subjected to the constraint $\dot{q}^{i}-(\rho_{\mathcal{D}})_{A}^{i}(q)y^{A}$, where we are replacing the control input in the cost function by the equation \eqref{equationlocal} that describes locally the solution of the fully actuated nonholonomic problem.

To derive the equations of motion of this variational problem with
constraints we can use standard variational calculus by defining the
extended Lagrangian
$$\widetilde{\mathcal{L}}(q^{i},y^{A},\dot{y}^{A},\lambda_i)=\mathcal{L}(q^{i},y^{A},\dot{y}^{A})+\lambda_{i}(\dot{q}^{i}-(\rho_{\mathcal{D}})_{A}^{i}(q)y^{A}),$$
and therefore the equations of motion determining necessary conditions for optimal extremals in the optimal control problem are $$
\dot{\lambda}_{i}=\frac{\partial\mathcal{L}}{\partial q^{i}}-\lambda_{j}\frac{\partial(\rho_{\mathcal{D}})_{A}^{j}}{\partial q^{i}}y^{A},\quad
\frac{d}{dt}\left(\frac{\partial\mathcal{L}}{\partial\dot{y}^{A}}\right)=\frac{\partial\mathcal{L}}{\partial y^{A}}-(\rho_{\mathcal{D}})_{A}^{i}\lambda_{i},\quad \dot{q}^{i}=(\rho_{\mathcal{D}})_{A}^{i}y^{A}.
$$

\subsubsection{Extension to underactuated systems:}  Consider the class of controllable underactuated nonholonomic mechanical systems, that is, controlled mechanical systems where the number of control inputs is less the rank of the linear subbundle spanned by the input sections. The search of necessary conditions for optimal extremals for the class of underactuated controlled nonholonomic mechanical systems can be done by using the same ideas than the fully actuated case. Assuming that the control subbundle $\mathcal{D}_{(c)}\subset E$ satisfies $\mathcal{D}_{(c)}=\hbox{span}\{e_{a}\}$ where
$\mathcal{D}=\hbox{ span}\left\{e_{a},e_{\alpha}\right\}=\hbox{ span}\left\{e_{A}\right\}$ and
$e_{A}$ are sections of $\tau_{\mathcal{D}}$, $\hbox{rank }\mathcal{D}_{(c)}=k<m=\hbox{ rank }\mathcal{D}$, a \textit{solution of an underactuated controlled nonholonomic problem} is an
admissible curve $\gamma:I\subset\R\to\mathcal{D}$ such that
$$\nabla_{\gamma(t)}^{\mathcal{G}^{\mathcal{D}}}\gamma(t)+grad_{\mathcal{G}^{\mathcal{D}}}V(\tau_{\mathcal{D}}(\gamma(t)))=u^{a}(t)e_{a}(\tau_{\mathcal{D}}(\gamma(t)).$$

Denote by $\{e^{a},e^{\alpha}\}$ the dual basis of
$\{e_{a},e_{\alpha}\}.$ This basis induces local
coordinates $(q^{i},y^{a},y^{\alpha})$ on $\mathcal{D},$ that is, if
$e\in\mathcal{D}$ then
$e=y^{A}e_{A}=y^{a}e_{a}+y^{\alpha}e_{\alpha}.$ Therefore, an
admissible curve has a local representation
$\gamma(t)=(q^{i}(t),y^{a}(t),y^{\alpha}(t))$ and optimal extremals for the underactuated nonholonomic problem are
characterized by admissible curves satisfying 
\begin{align*}
\Big{\langle}\nabla_{\gamma(t)}^{\mathcal{G}^{\mathcal{D}}}\gamma(t)+grad_{\mathcal{G}^{\mathcal{D}}}V(\tau_{\mathcal{D}}(\gamma(t))),e^{a}(\tau_{\mathcal{D}}(\gamma(t)) \Big{\rangle}=&u^{a}(t)\\
\Big{\langle}\nabla_{\gamma(t)}^{\mathcal{G}^{\mathcal{D}}}\gamma(t)+grad_{\mathcal{G}^{\mathcal{D}}}V(\tau_{\mathcal{D}}(\gamma(t))),e^{\alpha}(\tau_{\mathcal{D}}(\gamma(t))\Big{\rangle}=&0.
\end{align*}
Locally, the last equations read 
$$\dot{q}^i=(\rho_{\mathcal{D}})^i_A y^A,\quad
\dot{y}^{c}+\Gamma_{AB}^{c}y^{A}y^{B}+(\mathcal{G}^{\mathcal{D}})^{cB}\rho_{B}^{i}\frac{\partial
V}{\partial q^{i}}=u^{c},\quad\dot{y}^{\alpha}+\Gamma_{AB}^{\gamma}y^{A}y^{B}+(\mathcal{G}^{\mathcal{D}})^{\alpha
B}\rho_{B}^{i}\frac{\partial V}{\partial q^{i}}=0,$$ with $1\leq c\leq k=\hbox{rank }\mathcal{D}_{(c)}$ and $k+1\leq \alpha\leq m=\hbox{rank }\mathcal{D}$. The last set of equations is interpreted as constraints, therefore
we can denote by $\mathcal{M}\subset\mathcal{D}^{(2)}$ the
submanifold of $\mathcal{D}^{(2)}$ determined by these constraints.

Given a cost function $C:\mathcal{D}\times U\to\R$ the optimal
control problem consists on finding an admissible curve
$\gamma:I\subset\R\to\mathcal{D}$ solving the previous equations,
given boundary conditions on $\mathcal{D}$, and extremizing the cost functional
$\displaystyle{\mathcal{J}(\gamma(t),u(t))=\int_{0}^{T}C(\gamma(t),u(t))dt}.$ Solving the underactuated nonholonomic optimal control problem is equivalent to solving a
constrained variational problem determined by $\mathcal{L}:\mathcal{D}^{(2)}\rightarrow\mathbb{R}$
given in the selected coordinates by
\begin{equation*}\label{costunderactuated}\mathcal{L}(q^{i},y^{A},\dot{y}^{a})=C\left(q^{i},
y^{A},\dot{y}^{c}+\Gamma_{AB}^{c}y^{A}y^{B}+(\mathcal{G}^{\mathcal{D}})^{cB}\rho_{B}^{i}\frac{\partial
V}{\partial q^{i}}\right),\end{equation*} and subjected to the 
constraints $$\Phi^{\alpha}(q^{i},y^{A},\dot{y}^{\alpha})=\dot{y}^{\alpha}+\Gamma_{AB}^{\gamma}y^{A}y^{B}+(\mathcal{G}^{\mathcal{D}})^{\alpha
B}\rho_{B}^{i}\frac{\partial V}{\partial q^{i}}=0,\hbox{ and }
\dot{q}^{i}-\rho_{A}^{i}y^{A}=0.$$

To derive the equations of motion of this constrained  variational
problem we can use standard
variational calculus by extending the Lagrangian to $T\mathcal{D}$ with the
Lagrange multipliers $\lambda_{i}$ and $\overline{\lambda}_{\gamma}$
as
$$\widetilde{\mathcal{L}}(q^{i},y^{A},\dot{y}^{A}):=\mathcal{L}(q^{i},y^{A},\dot{y}^{a})+\lambda_{i}(\dot{q}^{i}-\rho_{A}^{i}(q)y^{A})+\overline{\lambda}_{\alpha}\Phi^{\alpha}(q^{i},y^{A},\dot{y}^{\alpha})$$ or by restricting the Lagrangian $\mathcal{L}$ to the submanifold $\mathcal{M}$. Proceeding with the first approach, we obtain that the necessary conditions for existence of extremals in the optimal control problem are determined by admissible curves
satisfying
\begin{align*}
0&=\dot{\lambda}_{i}+\lambda_{j}\frac{\partial\rho_{A}^{j}}{\partial q^{i}}y^{A}-\frac{\partial\mathcal{L}}{\partial q^{i}}-\overline{\lambda}_{\alpha}\frac{\partial \Phi^{\alpha}}{\partial q^{i}},\\
0&=\frac{d}{dt}\left(\frac{\partial\mathcal{L}}{\partial\dot{y}^{a}}\right)+\rho_{a}^{i}\lambda_{i}-\frac{\partial\mathcal{L}}{\partial y^{a}}-\overline{\lambda}_{\alpha}\left(\Gamma_{aB}^{\alpha}+\Gamma_{Ba}^{\alpha}\right)y^{B},\quad 0=\dot{q}^{i}-\rho_{A}^{i}y^{A}\\
0&=\dot{\overline{\lambda}}_{\alpha}+\frac{d}{dt}\left(\frac{\partial\mathcal{L}}{\partial\dot{y}^{\alpha}}\right)+\lambda_{i}\rho_{\alpha}^{i}-\frac{\partial\mathcal{L}}{\partial y^{\alpha}}-\overline{\lambda}_{\beta}\left(\Gamma_{\alpha B}^{\beta}+\Gamma_{B\alpha}^{\beta}\right)y^{B},\quad 0=\Phi^{\alpha}(q^{i},y^{A},\dot{y}^{A}).
\end{align*}

\subsection{Hamiltonian formulation of the optimal control problem}
 Necessary conditions of existence of extremals in the previous optimal control problem can be studied as a Hamiltonian problem on $T^{*}\mathcal{D}$ by defining the corresponding momenta for a constrained
(vakonomic) system (see for instance, in \cite{arnold2}, Section 4.2). The momenta are locally expressed as $$
p_i=\frac{\partial\widetilde{\mathcal{L}}}{\partial\dot{q}^{i}}=\frac{\partial\mathcal L}{\partial \dot{q}^i}+\lambda^j \frac{\partial  f_j}{\partial \dot{q}^i},\quad
p_A=\frac{\partial\widetilde{\mathcal{L}}}{\partial \dot{y}^{A}}=\frac{\partial\mathcal L}{\partial \dot{y}^A}+\lambda^j \frac{\partial  f^j}{\partial \dot{y}^A}$$
where $\widetilde{\mathcal L}$ is an arbitrary extension of ${\mathcal L}$ to $T{\mathcal D}$ using the constraints $f^j=\dot{q}^j-(\rho_{\mathcal{D}})^j_Ay^A=0$. If the map $\Psi:\mathcal{D}^{(2)}\times\mathbb{R}^{m}\to T^{*}\mathcal{D}$ locally given by $$\Psi(q^{i},y^{A},\dot{y}^{A},\lambda_{i})=(q^{i},y^{A},p_{i},p_{A})$$ is a local diffeomorphism, or equivalently, the matrix \[M=\left( \begin{array}{cc}
\frac{\partial^{2}\mathcal{L}}{\partial\dot{y}^{A}\partial \dot{y}^{A}} &\frac{\partial f^{j}}{\partial \dot{q}^{i}} \\
\left(\frac{\partial f^{j}}{\partial \dot{q}^{i}}\right)^{T} & 0 \end{array} \right)\]is non singular, the condition for local solvability of the constrained system is fulfilled and by the implicit function theorem one can locally define the Hamiltonian $\mathcal{H}:T^{*}\mathcal{D}\to\mathbb{R}$ as  $$\mathcal{H}(q^{i},y^{A},p_{i},p_{A})=p_{A}\dot{y}^{A}(q^{i},y^{A},p_{A}))+p_{i}\rho_{A}^{i}y^{A}-\mathcal{L}(q^{i},y^{A},\dot{y}^{A}(q^{i},y^{A},p_{A})).$$ Therefore, necessary conditions of existence of extremals  for the nonholonomic
optimal control problem are determined by the Hamiltonian system $(T^{*}\mathcal{D},\omega_{\mathcal{D}},\mathcal{H})$ where
$\omega_{\mathcal{D}}$ is the standard symplectic $2$-form on
$T^{*}\mathcal{D}.$ That is, these are determined by the equations
\begin{equation}\label{dynamic}
i_{X_{\mathcal{H}}}\omega_{\mathcal{D}}=d\mathcal{H}.
\end{equation}  Integral curves of $X_{\mathcal{H}}$ satisfies Hamilton's equations on $T^{*}\mathcal{D}$

$$
\dot{q}^{i}=\frac{\partial\mathcal{H}}{\partial p_{i}},\quad\quad \dot{y}^{A}=\frac{\partial\mathcal{H}}{\partial p_{A}},\quad
\dot{p}_{i}=-\frac{\partial\mathcal{H}}{\partial q^{i}},\qquad\dot{p}_{A}=-\frac{\partial\mathcal{H}}{\partial y^{A}}.$$
 
\begin{remark}
The previous equations specifying  the dynamics of the optimal control problem are the same in both frameworks, Lagrangian and Hamiltonian by using the identification provided by the momentum equations, and usually equations are given by a nonlinear system of equations, difficult to solve explicit. In order to integrate the equations, one approach is to apply numerical integrators for ordinary differential equations.
 
A numerical one-step method $y_{n+1}=\Psi_{h}(y_n)$ is called symplectic if, when applied to a Hamiltonian system on a symplectic space, the discrete flow $y\mapsto\Psi_{h}(y)$ is a symplectic transformation for all sufficiently small step sizes. We have seen that necessary conditions for the existence of optimal solutions in the optimal control problem can be seen as solutions of a Hamiltonian system on $T^{*}\mathcal{D}$ by defining a suitable Hamiltonian function $\mathcal{H}:T^{*}\mathcal{D}\to\mathbb{R}$, where solutions are  determined by the Hamiltonian vector field $X_{\mathcal{H}}$ for $\mathcal{H}$. 

One can then use standard methods for symplectic integration, such as symplectic Runge-Kutta methods, collocation methods, St$\ddot{\hbox{o}}$rmer-Verlet, Rattle and symplectic Euler methods  (see e.g. \cite{Hair}).  
For instance if we apply the St$\ddot{\hbox{o}}$mer-Verlet method for $\mathcal{H}:T^{*}\mathcal{D}\to\mathbb{R}$ we arrive to the algebraic system of equations
\begin{align*}
p_{n+1/2}&=p_{n}-\frac{h}{2}\nabla_{q}\mathcal{H}(q_n,p_{n+1/2}),\\
q_{n+1}&=q_n+\frac{h}{2}(\nabla_p \mathcal{H}(q_n,p_{n+1/2})+\nabla_p\mathcal{H}(q_{n+1},p_{n+1/2}),\\
p_{n+1}&=p_{n+1/2}-\frac{h}{2}\nabla_{q}\mathcal{H}(q_{n+1},p_{n+1/2}),\\
\tilde{p}_{n+1/2}&=\tilde{p}_{n}-\frac{h}{2}\nabla_{q}\mathcal{H}(q_n, p_n,y_n,\tilde{p}_{n+1/2}).
\end{align*} It would be interest to study the construction of symplectic integrators and compare the solutions with the ones obtained by applying variational integrators \cite{CoBl}. \hfill$\diamond$
\end{remark}
\begin{example}[Optimal control of the Chaplygin sleigh]
 A typical example of  nonholonomic  system  on  a  Lie  algebra,  we  study  the  Chaplygin  sleigh. The configuration space before reduction is the Lie group $G=SE(2)$ of the Euclidean motions of the 2-dimensional plane $\R^2$. It is well know that reduction of the system (see \cite{Bl}) gives rise to a Lagrangian defined on the Lie algebra of SE(2), denoted by $\mathfrak{se}(2)$. Elements of $\mathfrak{se}(2)$ are matrices of the form \[ \xi=\left(
\begin{array}{ccc} 0&\xi_3&\xi_1\\ -\xi_3&0&\xi_2\\ 0&0&0 \end{array} \right) \] and a basis of the Lie algebra $\mathfrak{se}(2)\cong \R ^3$ is
given by \[ E_1=\left( \begin{array}{ccc} 0&0&1\\ 0&0&0\\ 0&0&0 \end{array} \right) ,\qquad E_2=\left( \begin{array}{ccc} 0&0&0\\ 0&0&1\\ 0&0&0
\end{array} \right) ,\quad E_3=\left( \begin{array}{ccc} 0&-1&0\\ 1&0&0\\ 0&0&0 \end{array} \right). \] It is easy to check that $[E_3,E_1]=-E_2,
[E_2, E_3]=E_1, [E_1, E_2]=0$. An element $\xi \in \mathfrak{se}(2)$ is of the form  $\xi =v_1\, E_1+v_2\, E_2+\omega \, E_3$. The dynamics of the Chaplygin sleigh is described by the Euler-Poincar\'e-Suslov equations on $\mathfrak{se}(2)$. The Lagrangian function $L:\mathfrak{se}(2)\to\mathbb{R}$ is given by \[L(v_1, v_2, \omega)=\frac{1}{2}\left[ (J+m(a^2+b^2))\omega^2 + mv_1^2+m v_2^2-2bm\omega v_1-2am\omega v_2\right] \] where  $m$ and $J$ denotes the mass and moment of inertia of the sleigh relative to the contact point and $(a, b)$ represents the position
of the center of mass with respect to the body frame determined placing the origin at the contact point and the first coordinate axis in the
direction of the knife axis. 


The   system is subjected to the nonholonomic constraint determined by the linear subspace of $\mathfrak{se}(2)$: \[ {\mathcal D}=\{(v_1, v_2, \omega)\in se(2)\; |\; v_2=0\}\, . \] Instead of $\{E_1, E_2, E_3\}$ we take the basis of
$\mathfrak{se}(2)$ adapted to the decomposition ${\mathcal D}\oplus
{\mathcal D}^\perp$;  \[ \{X=E_3, Y=E_1, Z= -ma E_3-mab E_1+(J+ma^2) E_2\} \] ${\mathcal D}=\hbox{span }\{ X, Y\}$ and ${\mathcal D}^\perp=\hbox{span }\{ Z\}$. In the induced coordinates $(y_1, y_2)$ on ${\mathcal D}$ the restricted Lagrangian is \[ \ell(y_1, y_2)=\frac{1}{2}\left[ (J+m(a^2+b^2))(y_1^2 +
my_2^2-2bmy_1y_2\right]\; , \] and, given that $\lcf X, Y\rcf_{{{\mathcal D}}}=\frac{ma}{J+ma^2} X+\frac{mab}{J+ma^2}Y,\;$ it follows that ${\mathcal
C}^1_{12}=\frac{ma}{J+ma^2}$ and ${\mathcal C}^2_{12}=\frac{mab}{J+ma^2}$.

Using the relation between constant structures and Christofell symbols \eqref{relation}, the equations of motion are given by the Euler-Poincar\'e-Suslov equations on $\mathfrak{se}(2)$


$$ \dot{y}_1 =
\frac{mab}{J+ma^2}y_1^{2}-\frac{ma}{J+ma^2}y_1y_{2},\quad \dot{y}^2= \frac{ma(J+m(a^2+b^2))}{J+ma^2}y_1^2-\frac{mab}{J+ma^2}y_1y_2. $$

Next, consider an optimal control problem for the Chaplygin sleigh, where the control inputs are denoted by $u_1$ and $u_2$. The first one
corresponds to a force applied perpendicular to the center of mass
of the sleigh and the second one is corresponds with a force to control the heading direction. The controlled Euler-Poincar\'e-Suslov equations are 

\begin{align*} \dot{y}_1 -
\frac{mab}{J+ma^2}y_1^{2}+\frac{ma}{J+ma^2}y_1y_{2}&=u_{1},\\ \dot{y}_2- \frac{ma(J+m(a^2+b^2))}{J+ma^2}y_1^2+\frac{mab}{J+ma^2}y_1y_2&=u_2. \end{align*}
We want to find an admissible curve and control inputs satisfying the previous equations given boundary conditions on $\mathcal{D}$ and extremizing the cost functional $$\mathcal{J}=\int_{0}^{T}C(y_1,y_2,u_{1},u_{2})dt=
\frac{1}{2}\int_{0}^{T} (u_1^{2}+u_{2}^{2})dt.$$ Define the function $\mathcal{L}:\mathcal{D}^{(2)}\to\mathbb{R}$ in the induced coordinates $(y_1,y_2,\dot{y}_{1},\dot{y}_{2})$ by $\mathcal{L}(y_1,y_2,\dot{y}_{1},\dot{y}_2):=C(y_1,y_2,\dot{y}_{1},\dot{y}_{2})$, where  \begin{align*}
C(y_1,y_2,\dot{y}_{1},\dot{y}_{2})=&\frac{1}{2}\left(\dot{y}_1-\frac{mab}{J+ma^2}y_1^{2}+\frac{ma}{J+ma^2}y_1y_{2}\right)^{2}\\ 
&+\frac{1}{2}\left(\dot{y}_2- \frac{ma(J+m(a^2+b^2))}{J+ma^2}y_1^2+\frac{mab}{J+ma^2}y_1y_2\right)^{2}.\end{align*}

Denoting $K=J+m(a^2+b^2)$, necessary conditions of optimal curves on $\mathcal{D}$ are determined by the solutions of the nonlinear second-order system of equations

\begin{align*}
\frac{J+ma^2}{ma}\ddot{y}_{1}=&2by_1\dot{y}_{1}-\dot{y}_{1}y_2-\dot{y}_{2}y_1+(y_2-2by_1)\left(\dot{y}_1 -
\frac{mab}{J+ma^2}y_1^{2}+\frac{ma}{J+ma^2}y_1y_{2}\right)\\
&+\frac{ma}{J+ma^2}(by_2-2y_1K)\left(\dot{y}_{2}-\frac{maK}{J+ma^2}y_1^2+\frac{mab}{J+ma^2}y_1y_2\right),\\
\frac{J+ma^2}{ma}\ddot{y}_{2}=&y_1\dot{y}_1- \frac{ma(J+mab)y_1}{J+ma^2}y_1^2+\frac{may_1}{J+ma^2}y_1y_2\\
&+y_1b\left(\dot{y_2}-\frac{maK}{J+ma^2}y_1^2+\frac{mab}{J+ma^2}y_1y_2\right)+y_1\dot{y}_1K+b(\dot{y}_{1}y_2+y_1\dot{y}_{2}).
\end{align*}

It is strightfordward to see that the regularity condition ($\det M\neq 0$) is fulfilled. Therefore, by defining the corresponding momenta $$ p_1=\dot{y}_1-\frac{mab}{J+ma^2}y_1^2+\frac{ma}{J+ma^2}y_1y_2, \quad p_2=\dot{y}_2-\frac{maK}{J+ma^2}y_1^2+\frac{mab}{J+ma^2}y_1y_2$$ and denoting $(y_1,y_2,p_1,p_2)$ coordinates on $T^{*}\mathcal{D}$ we can define the Hamiltonian function describing the dynamics of the optimal control problem $$\mathcal{H}(y_1,y_2,p_1,p_2)=\frac{p_1^2}{2}+\frac{p_2^2}{2}+\frac{map_1}{J+ma^2}(by_1^2-y_1y_2)+\frac{map_2}{J+ma^2}(Ky_1^2+by_1y_2).$$ The corresponding Hamiltonian equations are
\begin{align*}
\dot{y}_{1}=&p_1+\frac{ma}{J+ma^2}(by_1^2-y_1y_2),\\
\dot{y}_{2}=&p_2+\frac{ma}{J+ma^2}(Ky_1^2+by_1y_2),\\
\dot{p}_{1}=&-\frac{map_1}{J+ma^2}(2by_1-y_2)-\frac{map_2}{J+ma^2}(2Ky_1+by_2),\\
\dot{p}_{2}=&\frac{map_1}{J+ma^2}y_1-\frac{map_2}{J+ma^2}by_1.
\end{align*} \hfill$\diamond$
\end{example}

\begin{example}[Optimal control of Euler-Poincar\'e-Suslov equations on $\mathfrak{so}(3)$]
Next we study optimal control for the Euler-Poincar\'e-Suslov equations on $\mathfrak{so}(3)$. Introducing controls in our picture, denoted by $u_1$ and $u_2$, the controlled dynamics is given by

$$ \dot{y}_1+\left(\frac{I_{13}}{I_{11}}y_1+\frac{I_{23}}{I_{22}}y_2\right)y_2=u_1,\quad  \dot{y}_2-\left(\frac{I_{13}}{I_{11}}y_1+\frac{I_{23}}{I_{22}}y_2\right)y_1=u_2. $$
We want to find an admissible curve and control inputs satisfying the previous equations given boundary conditions on $\mathcal{D}$ and extremizing the cost functional $$\mathcal{J}=\int_{0}^{T}C(y_1,y_2,u_{1},u_{2})dt=
\frac{1}{2}\int_{0}^{T} (u_1^{2}+u_{2}^{2})dt.$$ Define the function $\mathcal{L}:\mathcal{D}^{(2)}\to\mathbb{R}$ in the induced coordinates $(y_1,y_2,\dot{y}_{1},\dot{y}_{2})$ by the cost function $C:\mathcal{D}^{(2)}\to\mathbb{R}$ where 

$$
C(y_1,y_2,\dot{y}_{1},\dot{y}_{2})=\frac{1}{2}\left( \dot{y}_1+\left(\frac{I_{13}}{I_{11}}y_1+\frac{I_{23}}{I_{22}}y_2\right)y_2\right)^{2}
+\frac{1}{2}\left(\dot{y}_2-\left(\frac{I_{13}}{I_{11}}y_1+\frac{I_{23}}{I_{22}}y_2\right)y_1\right)^{2}.$$

Necessary conditions of optimal curves on $\mathcal{D}$ are determined by the solutions of the nonlinear second-order system of equations

\begin{align*}
\ddot{y}_{1}=&\frac{I_{13}y_2}{I_{11}}\dot{y}_1 +\frac{I_{13}}{I_{11}}y_2^2\left(
\frac{I_{13}}{11}y_1+\frac{I_{23}}{I_{22}}y_2\right)-y_2\left(\frac{I_{13}}{I_{11}}\dot{y}_{1}+\frac{I_{23}}{I_{22}}\dot{y}_{2}\right)-\dot{y}_{2}\left(\frac{I_{13}}{I_{11}}y_1+\frac{I_{23}}{I_{22}}y_2\right)\\
&+\left(\frac{I_{13}}{I_{11}}y_{1}^{2}+\frac{I_{23}}{I_{22}}y_1y_2-\dot{y}_{2}\right)\left(\frac{2I_{13}}{I_{11}}y_1+\frac{I_{23}}{I_{22}}y_2\right)\\
\ddot{y}_{2}=&\left(\dot{y}_1+\frac{I_{13}}{I_{11}}y_1y_2+\frac{I_{23}}{I_{22}}y_2^2 \right)\left(\frac{I_{13}}{I_{11}}y_1+\frac{2I_{23}}{I_{22}}y_2\right)+\dot{y}_1\left(
\frac{I_{13}}{11}y_1+\frac{I_{23}}{I_{22}}y_2\right)\\&+y_1\left(\frac{I_{13}}{I_{11}}\dot{y}_{1}+\frac{I_{23}}{I_{22}}\dot{y}_{2}\right)
-\frac{I_{23}y_1}{I_{22}}\left(\dot{y}_{2}-\frac{I_{13}}{I_{11}}y_{1}^{2}-\frac{I_{23}}{I_{22}}y_2y_1\right)
\end{align*}

The Hamiltonian description for the dynamics of the optimal control for the Euler-Poincar\'e-Suslov equations on $\mathfrak{so}(3)$ can be determined by defining the momenta $$p_1=\dot{y}_{1}+y_2\left(\frac{I_{13}}{I_{11}}y_1+\frac{I_{23}}{I_{22}}y_2\right),\quad p_2=\dot{y}_{2}-y_1\left(\frac{I_{13}}{I_{11}}y_1+\frac{I_{23}}{I_{22}}y_2\right)$$ since the regularity condition holds trivially.

Define the Hamiltonian function $\mathcal{H}:T^{*}\mathcal{D}\to\mathbb{R}$ in the induced coordinates $(y_1,y_2,p_1,p_2)$ by $$\mathcal{H}(y_1,y_2,p_1,p_2)=\frac{p_1^2}{2}+\frac{p_2^2}{2}+p_2y_1\left(\frac{I_{13}}{I_{11}}y_1+\frac{I_{23}}{I_{22}}y_2\right)-p_1y_2\left(\frac{I_{13}}{I_{11}}y_1+\frac{I_{23}}{I_{22}}y_2\right).$$

The Hamiltonian equations describing the necessary conditions for optimal trajectories in the optimal control problem are
\begin{align*}
\dot{y}_{1}=&p_1-y_2\left(\frac{I_{13}}{I_{11}}y_1+\frac{I_{23}}{I_{22}}y_2\right),\\
\dot{y}_{2}=&p_2+y_1\left(\frac{I_{13}}{I_{11}}y_1+\frac{I_{23}}{I_{22}}y_2\right),\\
\dot{p}_{1}=&p_1y_2\frac{I_{13}}{I_{11}}-\frac{I_{13}p_2}{I_{22}}\left(\frac{I_{13}}{I_{11}}y_1+\frac{I_{23}}{I_{22}}y_2\right),\\
\dot{p}_{2}=&p_1\left(\frac{I_{13}}{I_{11}}y_1+\frac{I_{23}}{I_{22}}y_2\right)+p_2\left(\frac{I_{13}}{I_{11}}y_1+\frac{I_{23}}{I_{22}}y_2\right)+\frac{I_{23}p_1y_2}{I_{22}}\left(\frac{I_{13}}{I_{11}}y_1+\frac{I_{23}}{I_{22}}y_2\right)-\frac{p_2y_1I_{23}}{I_{22}}.
\end{align*}

 \hfill$\diamond$
\end{example}

\end{document}